\documentclass{amsart}
\usepackage{graphicx}
\usepackage{amsfonts, amsmath, amsfonts, amssymb}
\newtheorem{thm}{Theorem}[section]
\newtheorem{cor}[thm]{Corollary}
\newtheorem{lem}[thm]{Lemma}
\newtheorem{prop}[thm]{Proposition}
\theoremstyle{definition}
\newtheorem{defn}[thm]{Definition}
\theoremstyle{remark}
\newtheorem{rem}[thm]{Remark}
\numberwithin{equation}{section}
\newenvironment{prf}{\noindent{\bf Proof}}{\\ \hspace*{\fill}$\Box$ \par}
\newenvironment{skprf}{\noindent{\\ \bf Sketch of Proof}}{\\ \hspace*{\fill}$\Box$ \par}

\begin{document}

\title[Note on a theorem of Munkres]{A note on a theorem of Munkres}%
\author{Emil Saucan}%
\address{Department of Mathematics, Technion, Haifa and Software Engineering Department, Ort Braude College, Karmiel, Israel}%
\email{semil@tx.technion.ac.il}%

\thanks{This paper represents an outshoot of the authors Ph.D. Thesis written under the supervision of Prof. Uri Srebro}%

\date{14.3.2004.}%
\begin{abstract}
We prove that given a $\mathcal{C^\infty}$ Riemannian manifold with boundary, any fat triangulation of the
boundary can be extended to the whole manifold. We also show that this result holds extends to $\mathcal{C}^1$
manifolds, and that in dimensions $2,3$ and $4$ it also holds for $PL$ manifolds. We employ the main result to
prove that given any orientable $\mathcal{C^\infty}$ Riemannian manifold  with boundary admits quasimeromorphic
mappings onto $\widehat{\mathbb{R}^n}$. In addition some generalizations are given.
\end{abstract}
\maketitle
\section{Introduction}
 The existence of triangulations for $\mathcal{C}^1$  manifolds without boundary is known since the classical work of Whithead (\cite{wh},
 1940).
\\ This result was  extended in 1960 by  Munkres (\cite{mun}) to include $\mathcal{C}^{r},\; 1\leq r \leq \infty$, manifolds with
boundary. To be more precise, he proved that any $\mathcal{C}^{r}$ triangulation of the boundary can be extended
to a $\mathcal{C}^{r}$ triangulation of the whole manifold.
\\ Yet earlier, in 1934-1935 (even before Whithead's work) Cairns\footnote{\, Although far better known and widely cited, Whitehead's work is rooted in
Cairns' studies, to whom it gives due credit in the very opening phrase: "{\tt This paper is suplementary to S.S.
Cairns' work on the triangulation ... of manifolds of class $\mathcal{C}^1$.}"}  (\cite{ca1}, \cite{ca2}) proved
the existence of fat triangulations for compact $\mathcal{C}^1$ manifolds and for manifolds with boundary having a
finite number of compact boundary components. Moreover, his triangulations were {\it fat} (see definition bellow),
something which Munkres' method achieved only away from the boundary (see Section 2.2).
\\ Unfortunately, it seems that little interest existed during the following decades, for studying generalizations of the results
above (\cite{fe} representing a notable exception).
\\ The interest in the existence of fat triangulation was rekindled by the study of quasiregular and
quasimeromorphic functions, since the existence of fat triangulations is crucial in the proof of existence of
guasiregular (quasimeromorphic) mappings (see \cite{ms2}, \cite{tu}) and in 1992  Peltonen (\cite{pe}) proved the
existence of fat triangulations for $\mathcal{C}^{\infty}$ Riemannian manifolds, using methods partially based
upon another technique of Cairns (originally developed for triangulating manifolds of class $\geq \mathcal{C}^2$).
\\ In this paper we 
extend Munkres' Theorem to the case of fat triangulations of  manifolds (and orbifolds -- see Section 5.2.) with
or without boundary and we show how to apply this main result in order to prove the existence of quasimeromorphic
functions. Our main result is the following Theorem:

\begin{thm}
Let $M^{n}$ be an $n$-dimensional $\mathcal{C}^{\infty}$ manifold with boundary. Then any fat triangulation of
$\partial M^{n}$ can be extended to a fat triangulation of $M^{n}$.
\end{thm}

\hspace*{-0.35cm}where a fat triangulations is defined as follows:

\begin{defn} Let $\tau \subset \mathbb{R}^n$ ; $0 \leq k \leq n$ be a $k$-dimensional simplex.
The {\it fatness}  $\varphi$ of $\tau$ is defined as being:
\begin{equation}
\varphi = \varphi(\tau) = \inf_{\hspace{0.4cm}\sigma < \tau
\raisebox{-0.25cm}{\hspace{-0.9cm}\mbox{\scriptsize$dim\,\sigma=l$}}}\!\!\frac{Vol(\sigma)}{diam^{l}\,\sigma}
\end{equation}
The infimum is taken over all the faces of $\tau$, $\sigma < \tau$, and $Vol_{eucl}(\sigma)$ and $diam\,\sigma$
stand for the Euclidian $l$-volume and the diameter of $\sigma$ respectively. (If $dim\,\sigma = 0$, then
$Vol_{eucl}(\sigma) = 1$, by convention.)
\\ A simplex $\tau$ is $\varphi_0${\it -fat}, for some $\varphi_0 > 0$, if $\varphi(\tau) \geq \varphi_0$. A triangulation (of a submanifold of $\mathbb{R}^n$) $\mathcal{T} = \{ \sigma_i \}_{i\in \bf I}$ is
$\varphi_0${\it -fat} if all its simplices are $\varphi_0$-fat. A triangulation $\mathcal{T} = \{ \sigma_i
\}_{i\in \bf I }$ is {\it fat} if there exists $\varphi_0 \geq 0$ s.t. all its simplices are $\varphi_0${\it-fat};
$\forall i \in \bf{I}. $
\end{defn}

\begin{rem} There exists a constant $c(k)$ that depends solely upon the dimension $k$ of  $\tau$ s.t.
\begin{equation}
\frac{1}{c(k)}\cdot\varphi(\tau) \leq \min_{\hspace{0.4cm}\sigma < \tau
\raisebox{-0.25cm}{\hspace{-0.9cm}\mbox{\scriptsize$dim\,\sigma=l$}}}\hspace{-0.3cm}\measuredangle(\tau,\sigma)
\leq c(k)\cdot\varphi(\tau)\,,
\end{equation}
and
\begin{equation}
\varphi(\tau) \leq \frac{Vol(\sigma)}{diam^{l}\,\sigma} \leq c(k)\cdot\varphi(\tau)\,;
\end{equation}
where $\measuredangle(\tau,\sigma)$ denotes the  ({\it internal}) {\it dihedral angle} of $\sigma < \tau$. (For a
formal definition, see \cite{cms}, pp. 411-412, \cite{som}.)
\end{rem}

\begin{rem}
The definition above is the one introduced in \cite{cms}. For different, yet equivalent definitions of fatness,
see \cite{ca1}, \cite{ca2}, \cite{pe}, \cite{tu}.
\end{rem}

The idea of the proof of Theorem 1.1. is first to build two fat triangulations: $\mathcal{T}_{1}$ of a product
neighbourhood $N$ of $\partial M^n$ in $int\,M^n$ and $\mathcal{T}_{2}$ of $int\, M^n$, and then to "mash" the two
triangulations into a new triangulation $\mathcal{T}$, while retaining their fatness.

While the mashing procedure of the two triangulations is basically that developed in the original proof of
Munkres' theorem, the triangulation of $\mathcal{T}_{1}$ was modified, in order to ensure the fatness of the
simplices of $\mathcal{T}_{1}$. The existence of the second triangulation is assured by Peltonen's result.
\\ Thus our main efforts are dedicated to the task of fattening the newly obtained triangulation into a new
fat triangulation. The technique we employ is essentially the one developed  in \cite{cms}.\footnote{\, For a more
direct approach in dimensions 2 and 3 see \cite{s1}. Also, for the treatment of the same problem in the context of
Computational Geometry, see  \cite{e}.}
\\ Once a fat triangulation of an orientable manifold $ M^n$ is provided, the construction of the required quasimeromorphic mapping is
canonical (see \cite{al}, \cite{pe}, \cite{ms2}, \cite{tu}) and is based upon the so called "Alexander Trick",
which we present here succinctly (in a nutshell): one starts by constructing a suitable triangulation of $M^n$.
Since $M^n$ is orientable, an orientation consistent with the given triangulation (i.e. such that two given
$n$-simplices having a $(n-1)$-dimensional face in common will have opposite orientations) can be chosen. Then one
quasiconformally maps the simplices of the triangulation into $\widehat{\mathbb{R}^n}$ in a chess-table manner:
the positively oriented ones onto the interior of the standard simplex in $\mathbb{R}^n$ and the negatively
oriented ones onto its exterior. If the dilatations of the quasiconformal maps constructed above are uniformly
bounded -- which condition is fulfilled if the simplices of the triangulation are of uniform fatness -- then the
resulting map will be quasimeromorphic.

 This paper is organized as follows: in Section 2 we present the main techniques we employ: Peltonen's method of triangulating $int\,M^n$  and the Proof of Munkres' Theorem on the extension of the triangulation of $\partial M^n$ to $int\, M^n$.
Section 3 is dedicated to the main task of fattening the common triangulation. In Section 4 we show how to apply
the main result in the construction of a quasimeromorphic mappings from $M^n$ to $\widehat{\mathbb{R}^n}$.
Finally, in Section 5 we propose some generalizations.

\section{Extending $\mathcal{T}_1$ to $int\, M^n$}

\subsection{Peltonen's Technique} Peltonen's method is an extension of one due to Cairns, developed in order to triangulate $\mathcal{C}^{2}$-compact manifolds
(\cite{ca3}). It is based on the subdivision of the given manifold into a closed cell complex generated by a
Dirichlet (Voronoy) type partition whose vertices are the points of a maximal set that satisfy a certain density
condition. We give below a sketch of the Peltonen's method, refering the interested reader to \cite{pe} for the
full details.
\\The construction devised by Peltonen consists of two parts:
\\{\em Part 1} This part proceeds in two steps:
\\  \hspace*{0.3cm}{\em Step A} Build an exhaustation $\{E_i\}$ of $M^n$, generated by the pair $(U_i,\eta_i)$, where:
\begin{enumerate}
\item $U_i$ is the relatively compact set $E_i \setminus E_{i-1}$ and
\item $\eta_i$ is a number that controls the fatness of the simplices of the triangulation of $E_i$, that will be constructed in Part 2, such that it will not differ to much
on adjacent simplices, i.e.:
\\ (i) The sequence $(\eta_i)_{i\geq1}$ descends to $0$\,;
\\ (ii) $2\eta_i \geq \eta_{i-1} \,.$
\end{enumerate}
 \hspace*{0.3cm}{\em Step B}
\begin{enumerate}
\item Produce a maximal set $A$, $|A| \leq \aleph_0$, s.t. $A \cap U_i$ satisfies:
\\ (i) a density condition, and
\\ (ii) a "gluing" condition (for $U_i, U_{i=1}$).
\item Prove that the Dirichlet complex $\{\bar{\gamma}_i\}$ defined by the sets $A_i$ is a cell complex and
every cell has a finite number of faces (so it can be triangulated in a standard manner).
\end{enumerate}
{\em Part 2} \,Consider first the dual complex $\Gamma$ and prove that it is a Euclidian simplicial complex with a
"good" density, then   project $\Gamma$ on $M^n$ (using the normal map). Finally, prove that the resulting complex
can be triangulated by fat simplices.

\begin{rem}
In the course of Peltonen's construction $M^n$ is presumed to be isometrically embedded in some
$\mathbb{R}^{N_1}$, where the existence of $N_1$ is guaranteed by Nash's Theorem (see \cite{pe}, \cite{spi5}).
\end{rem}

\subsection{The Extension of $\mathcal{T}_1$ to $int\, M^n$}
We first establish some notations and definitions:
\\ Let $K$ denote a simplicial complex, and let $K' < K$ denote a subcomplex of $K$.

\begin{defn}
Let $f_i:K_i \stackrel{\sim}{\rightarrow} \mathbb{R}^n,\; i=1,2$ be s.t. $f(|K_i|)$ is closed.
\\ We say that $(K_1,f_1), (K_2,f_2)$ {\em intersect in a subcomplex} iff:
\\ (i) $f_i^{-1}\big(f_1(|K_1|) \cap f_2(|K_2|)\big) = |L_i|$\,; where $L_i < K_i\,,\; i=1,2$.
\\ and
\\ (ii) $f_2^{-1}\circ f_1: L_1 \rightarrow L_2$ is a linear isomorphism.\footnote{\,i.e. (i) $f:|L_1| \stackrel{\sim}{\rightarrow} |L_2|$ and (ii) $f|_{\sigma}$ is linear, $\forall \sigma \in L_1$.}
\end{defn}

\begin{defn}
Let $L < K$. $L$ is called {\em full} iff $\sigma \cap L$ either is a face of $\sigma$ or else it is empty;
$\forall \sigma \in K$.
\end{defn}

\begin{rem}
$L$ is full $\Longleftrightarrow \partial \sigma \cap L \neq \partial \sigma;\; \forall \sigma \in K$. 
\end{rem}

If $(K_1,f_1), (K_2,f_2)$ intersect in a full subcomplex, then there exist a complex $K$ and a homeomorphism $f:K
\rightarrow \mathbb{R}^n$ s.t. the following diagram is comutative:  

\[ \hspace*{-2.6cm}K_1 \]
\begin{picture}(10,20)(1,2)
 \put(140,21){\vector(0,-1){31}} \put(153,27){\vector(3,-2){57}} \put(175,17){$f_1$} \put(130,5){$i_1$}
\end{picture}

\[ \hspace*{-0.1cm}K  \hspace{2.3cm} \mathbb{R}^n\]
\begin{picture}(10,20)(1,2)
 \put(145,27){\vector(1,0){63}} \put(175,31){$f$}
\end{picture}

\[ \hspace*{-2.6cm}K_2 \]
\begin{picture}(10,20)(1,2)
 \put(140,39){\vector(0,1){31}} \put(153,31){\vector(3,2){57}} \put(175,35){$f_2$} \put(130,50){$i_2$}
\end{picture}

Here $i_1, i_2$ are linear isomorphisms. The pair $(K,f)$ is unique up to isomorphism.

\begin{defn}
Let $(K_1,f_1), (K_2,f_2)$ and $(K,f)$ be as above. Then $(K,f)$ is called the {\em union} of $(K_1,f_1)$ and
$(K_2,f_2)$.
\end{defn}

\begin{defn}
Let $f:K \rightarrow \mathbb{R}^n$ be a $\mathcal{C}^r$ map, and let $\delta:K \rightarrow \mathbb{R}^*_+$ be a
continuous function. Then $g:|K| \rightarrow \mathbb{R}^n$ is called a $\delta${\it-approximation to} $f$ iff:
\\ (i) There exists a subdivision $K' < K$ s.t. $g \in \mathcal{C}^r(K',\mathbb{R}^n)$\,;
\\ (ii) $d_{eucl}\big(f(x),g(x)\big) < \delta(x)$\,, $\forall x \in |K|$\,;
\\ (iii) $d_{eucl}\big(df_a(x),dg_a(x)\big) \leq \delta(a)\cdot d_{eucl}$\,, $\forall a \in |K|, \forall x \in \overline{St}(a,K')$
\end{defn}

\begin{defn}
Let $K < K'$, $U = \raisebox{0.05cm}{\mbox{$\stackrel{\circ}{U}$}}$, and let $f \in \mathcal{C}^r(K,\mathbb{R}^n),
\; g  \in \mathcal{C}^r(K,\mathbb{R}^n)$.  g is called a $\delta${\it-approximation} of $f$ (on $U$) iff
conditions (i) and (ii) of Definition 2.1. hold for any $a \in U$.
\end{defn}

\begin{defn}
Let $K < K'$ and let $f \in \mathcal{C}^r(K,\mathbb{R}^n), \; g \in \mathcal{C}^r(K',\mathbb{R}^n)$ be
non-degenerate\footnote{\,i.e. $rank(f|_{\sigma}) = dim\,\sigma, \; \forall \sigma \in K$.} mappings and let $U =
\raisebox{0.05cm}{\mbox{$\stackrel{\circ}{U}$}} \subset |K|$. $g$ is called an $\alpha${\it-approximation} (of $f$
on $U$) iff:
\begin{equation}
\angle\big(df_a(x),dg_a(x)\big) \leq \alpha\,;\; \forall a \in U, \forall x \in \overline{St}(a,K'),\; a \neq x.
\end{equation}
\end{defn}

We now bring Munkres' Theorem. While we will initially apply it for $\mathcal{C}^\infty$ manifolds, we give the
proof for the general case of $\mathcal{C}^r$ manifolds, $1\leq r\leq \infty$. We modify the original construction
so the triangulation of a certain neighbourhood of $\partial M^n$ will be fat.

\begin{thm}[\cite{mun},  10.6]
Let $M^{n}$ be a $\mathcal{C}^{r}$-manifold with boundary. Then any $\mathcal{C}^{r}$-triangulation of $\partial
M^{n}$ can be extended to a $\mathcal{C}^{r}$-triangulation of $M^{n}$, $1\leq r\leq \infty$.
\end{thm}

\begin{proof}
Let $f:J \rightarrow \partial M^n$ be a $\varphi_{_{\partial M}}$-fat $\mathcal{C}^r$ triangulation, for some
$\varphi_{_{\partial M}}$. We construct a triangulation of $|J| \times [0,1)$ in the following way: \\ If $J$ is
isometrically embedded in $\mathbb{R}^{N_2}$, we consider (in $\mathbb{R}^{N_2}$) the cells of type:

\begin{equation}
\sigma_{1,n} = \sigma \times \Big[\frac{k}{n_0}\,,\,\frac{k+1}{n_0}\Big]\,; \; k = 1, \ldots ,n_0 - 1\,.
\end{equation}

and

\begin{equation}
\sigma_{2,n} = \sigma \times \Big\{\frac{k}{n_0}\Big\}\,; \forall \sigma \in J.
\end{equation}

Let $K$ denote the resulting simplicial complex: $|K| = |J| \times [0,1)$. The cells of the complex above may be
divided in simplices without subdividing the cells of type $\sigma_{2,n}$\,. (See Fig. 1 for the case $N_2 = 2$.)

\begin{figure}[h]
\begin{center}
\includegraphics[scale=0.3]{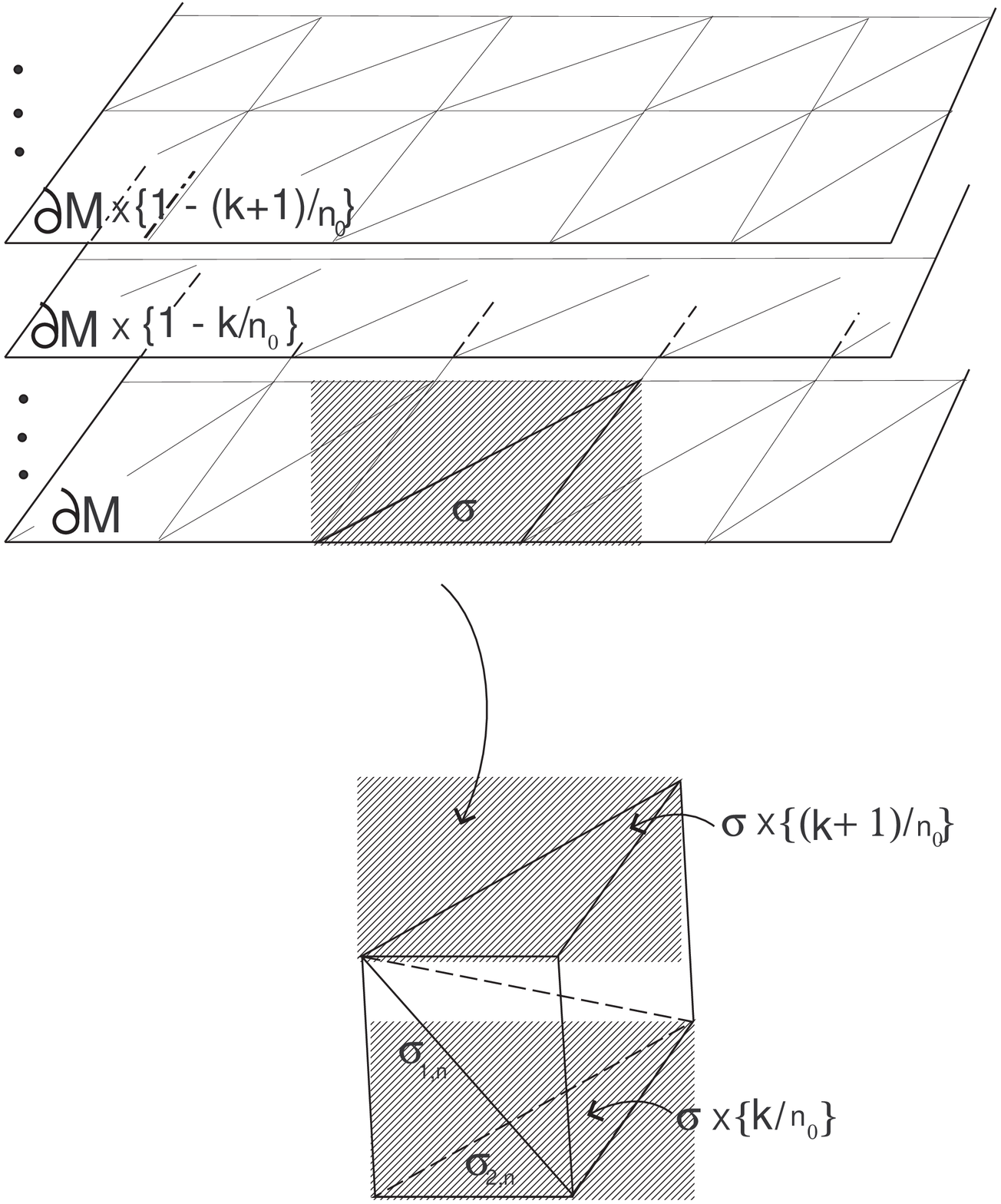}
\end{center}
\caption{ }
\end{figure}

For reasons that will become clear in the course of the proof of Theorem 3.6.\,, we choose $n_0$  such that the
fatness of any simplex $\sigma \in K$ is $\geq \varphi_0$, for some $\varphi_0$\footnote{\,$\varphi_0$ depends
upon $\varphi_{_{\partial M}}$ and $\varphi_{_{int\,M}}$\,.} and such that $diam\,\sigma \leq diam\,\tau ,\;
\forall \sigma \in K_0,\; \tau \in L_0$\footnote{\, To attain this inequalities, further subdivision may be
necessary -- their number depending upon the respective "$\eta_i$"-s given by Peltonen's construction.}\,, where
$K_0, L_0$ are defined as follows:

Let $K_0$ be the subcomplex of $K$ s.t. $|K_0| = |J| \times \big[0\,,\,\frac{k_1}{n_0}\big],\; k_1 =
[\frac{5n_0}{6}]$\,; and let
\\$\psi:\partial M^n \times [0,1) \rightarrow M^n$ be a product neighbourhood of $\partial M^n $ (in $M^n$). (Here $M^n$ is supposed to be embedded in $N = \max\{N_1,N_2\}$\,.)
Then, if $g$ makes the following diagram commutative:
\newpage
\[J \times [0,1) \hspace{1.7cm} \partial M^n \times [0,1]\]
\begin{picture}(10,20)(1,2)
\put(158,29){\vector(1,0){33}} \put(147,24){\vector(2,-3){23}} \put(210,24){\vector(-2,-3){23}} \put(150,5){$g$}
\put(203,5){$\psi$} \put(161,32){$f \times id$}
\end{picture}
\vspace*{0.3cm}\[M^n\] \vspace*{0.5cm} then $g|\raisebox{-0.3cm}{\scriptsize\mbox{$K_0$}}$ is a $\mathcal{C}^r$
embedding s.t.
\\ (i) $g(K_0) = \overline{g(K_0)}$ (in $M^n$)
\\ and
\\ (ii) $\psi\big(\partial M^n \times [0,\frac{k_2}{n_0})\big) \subset int\, g(K_0),\; k_2 = [\frac{4n_0}{5}]$.
\\ Now, if $h:L \rightarrow M^n$ is a $\mathcal{C}^r$ triangulation of $int\,M^n$, then, by further (eventual)
subdivision, we may suppose\footnote{\,see \cite{mun}.} that: $\sigma' \cap \,\psi\big(\partial\,M^n \times
[0,\frac{k_3}{n_0}]\big) = \emptyset,\; k_3 = [\frac{3n_0}{4}]\,; \forall\, \sigma' \in L$; \newline $\sigma' \cap
\psi\big(\partial\,M^n \times \{\frac{k_2}{n_0}\}\big) \neq \emptyset$.
\\ Let $L_0$ be the complex given by:
\begin{equation}
\left \{ \begin{array}{lll}
  L_0^i = \big\{\sigma \in L \,\big|\, h(\sigma) \cap \big(M^n
\,\backslash\,\psi\big(\partial\,M^n
\times [0,\frac{k_2}{n_0})\big)\big) \neq \emptyset \big\}\,; \\
  L_0^f = \big\{{\rm faces\; of }\; \sigma \,\big|\, \sigma \in L_0^i\big\}\,; \\
  L_0 = L_0^i \cap L_0^f\,.
      \end{array}
\right.
\end{equation}

Then, by \cite{mun}, Theorem 10.4, (see also Fig. 2\,\footnote{\, after \cite{mun}.}) $\exists\, g':K'_0
\rightarrow M^n,\; h': L'_0 \rightarrow M^n$; where $g'$ is a $\delta$-approximation of $g$ and $h$ is a
$\delta$-approximation of $h$, s.t.
\\ (i) $g'(K'_0) \cap h'(L'_0)$ is full
\\and
\\ (ii) The union of $(K'_0,g')$ and $(L'_0,h')$ is an embedding.

\begin{figure}[h]
\begin{center}
\includegraphics[scale=0.4]{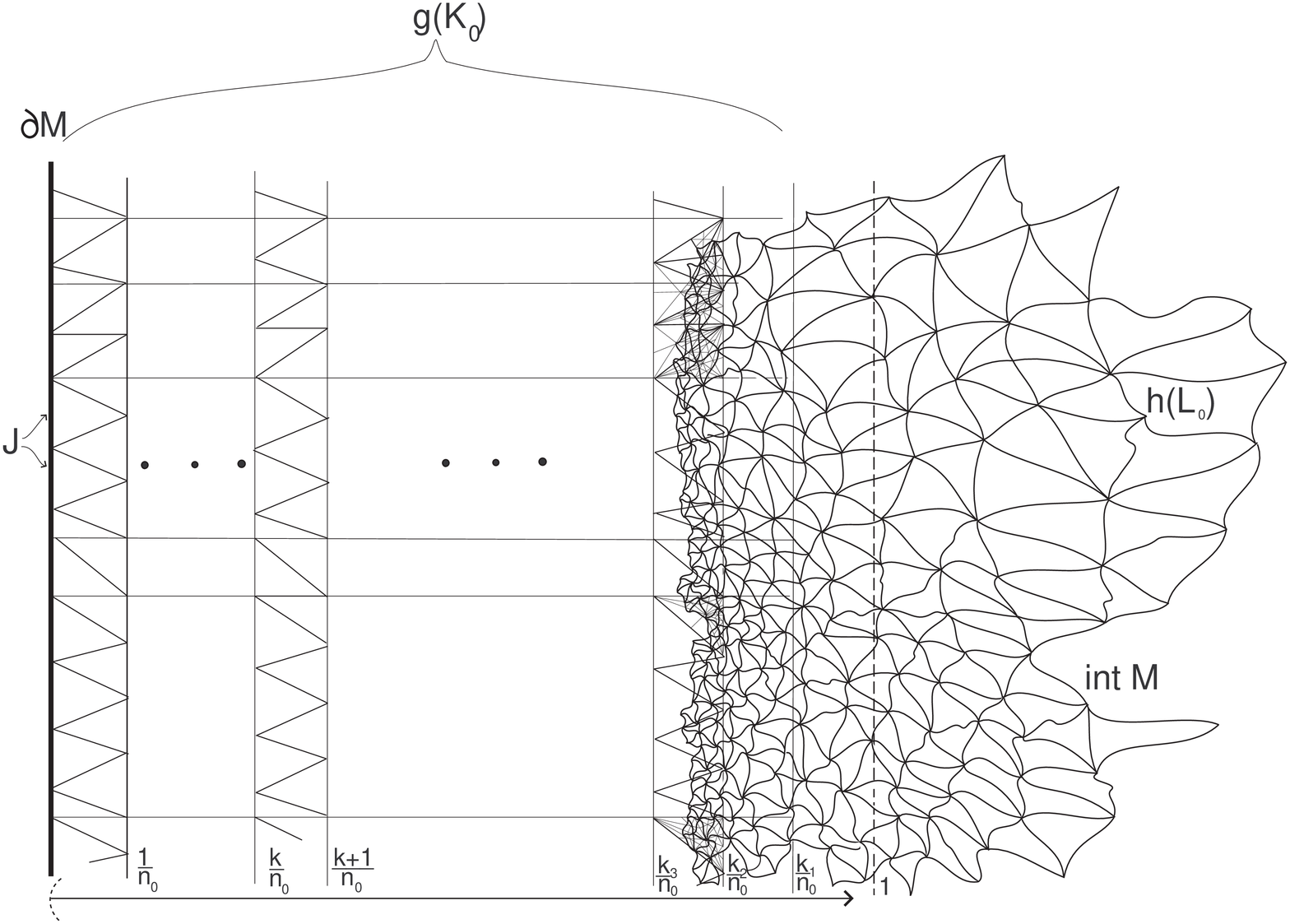}
\end{center}
\caption{ }
\end{figure}

Also, by applying again \cite{mun} Theorem 10.4, we may suppose that
\\ (a) $K'_0\big|\raisebox{-0.3cm}{\scriptsize\mbox{$|J|\times [0,\frac{k_4}{n_0}]$}} \equiv K_0\big|\raisebox{-0.3cm}{\scriptsize\mbox{$|J|\times [0,\frac{k_4}{n_0}]$}}$
\\ (b) $g'_0\big|\raisebox{-0.3cm}{\scriptsize\mbox{$|J|\times [0,\frac{k_4}{n_0}]$}} \equiv g_0\big|\raisebox{-0.3cm}{\scriptsize\mbox{$|J|\times [0,\frac{k_4}{n_0}]$}};\; k_4 = [\frac{n_0}{2}]$

Then $(K'_0,g') \cup (L'_0,h')$ will be the sought for triangulation, but only if the following condition also
holds:
\begin{equation}
g'\big(|K_0|\big) \cup h'\big(|L_0|\big) = M^n\,.
\end{equation}

But this condition also takes hold in our case, by virtue of a more general result about topological manifolds
(see \cite{mun}, pp. 36-38, 105).

\end{proof}

\section{Fattening Triangulations}

First let us establish some definitions and notations:

\begin{defn}
Let $\sigma_i \in K$, $dim\,\sigma_i = k_i$, $i=1,2$; s.t. $diam\,\sigma_1 \leq diam\,\sigma_2$. We say that
$\sigma_1, \sigma_2$ are $\delta${\em-transverse} iff
\\ (i) $dim(\sigma_1 \cap \sigma_2) = \max(0,k_1+k_2-n)$;
\\ (ii) $0 < \delta < \measuredangle(\sigma_1,\sigma_2)$;
\\ and if $\sigma_3 \subset \sigma_1,\,\sigma_4 \subset \sigma_2$, s.t. $dim\,\sigma_3 + dim\,\sigma_4 < n =
dim\,K$, then
\\ (iii) $dist(\sigma_3, \sigma_4) > \delta\cdot\eta_1$.
\\ In this case we write: $\sigma_1 \pitchfork_{\delta} \sigma_2$.
\end{defn}

We begin by triangulating and fattening the intersection of two individual simplices belonging to the two given
triangulations, respectively. Given two closed simplices $\bar{\sigma_1}, \bar\sigma_2$, their intersection (if
not empty) is a closed, convex polyhedral cell: $\bar\gamma = \bar\sigma_1 \cap \bar\sigma_2$. One canonically
triangulates  $\bar\gamma$ by using the {\em barycentric subdivison} $\bar\gamma^\ast$ of $\bar\gamma$, defined
inductively upon the dimension of the cells of $\partial\gamma$  in the following manner: for each cell $\beta
\subset \partial\gamma$, choose an interior point $p_{\beta} \in int\,\beta$ and construct the join
$J(p_\beta,\beta_i),\; \forall \beta_i \subset \partial \beta$.\footnote{\, If $dim\,\beta = 0$ or $dim\,\beta =
1$, then $\beta$ is already a simplex.}

We first show that if the simplices are fat and if they intersect $\delta$-transversally, then one can choose the
points s.t. the barycentric subdivision $\bar\gamma^\ast$ will be composed of fat simplices. More precisely, we
prove the following Proposition:

\begin{prop}
Let $\sigma_1, \sigma_2 \subset \mathbb{R}^m$, where $m = \max(dim\,\sigma_1,dim\,\sigma_2)$, s.t.
\\ $d_1 = diam\,\sigma_1 \leq d_2 = diam\,\sigma_2$, and s.t. $\sigma_1, \sigma_2$ have common fatness $\varphi_0$.
\\ If $\sigma_1 \pitchfork_{\delta} \sigma_2$, then there exists $c = c(m,\varphi_0,\delta)$

\begin{enumerate}
\item If $\sigma_3 \subset \bar\sigma_1,\,\sigma_4 \subset \bar\sigma_2$ and if $\sigma_1 \cap \sigma_2 \neq \emptyset$, then $\sigma_1 \cap \sigma_2 =
\gamma_0$ is an ($k_3+k_4-m$)-cell, where $k_3 = dim\,\sigma_3, k_4 = dim\,\sigma_4$ and:

\begin{equation}
Vol_{eucl}(\gamma_0) \geq c\cdot d_1^{k_3+k_4-m}\;.
\end{equation}

\item $\forall\, \gamma_0$ as above, $\exists\, p \in \gamma_0$, s.t.
\begin{equation}
dist(p,\partial\gamma_0) > c\cdot d_1\;.
\end{equation}

\item If the points employed in the construction of $\bar\gamma^\ast$ satisfy the condition {\rm (3.2)} above, then each
$l$-dimensional simplex $\tau \in \gamma^\ast$ satisfies the following
\\inequalities:

\begin{equation}
\varphi_l \geq Vol_{eucl}(\tau)/d_1^{\,l} \geq  c\cdot d_1\;.
\end{equation}

\end{enumerate}

\end{prop}

\begin{proof}

First, consider the following remarks:

\begin{rem}
 The following sets are compact:
\\ $S_1 = \{\sigma_1\,|\, diam\,\sigma_1 = 1\,, \varphi(\sigma_1) \geq \varphi_0\}$, $S_2 = \{\sigma_2\,|\, diam\,\sigma_2 = 2(1+\delta)\,, \varphi(\sigma_2) \geq
\varphi_0\}$,

$S(\phi_0,\delta) \subset S_1 \cap S_2\,, \; S(\phi_0,\delta) = \{(\sigma_1, \sigma_2)\,|\, \exists v_0, {\rm
s.t.} v_0 \in \sigma_1, \forall \sigma_1 \in S_1 \cap S_2\}$.
\end{rem}

\begin{rem}
There exists a constant $c(\varphi)$ s.t. $\mathcal{S} = \mathcal{S}'$, where
\\ $\mathcal{S} = \{\sigma_1 \cap \sigma_2\,|\, diam\,\sigma_2 \leq d_2\},\;  \mathcal{S}' = \{\sigma_1 \cap \sigma_2\,|\, diam\,
c(\varphi)(1+\delta)d_1\}$,
\\ i.e. the sets of all possible intersections remains unchanged under controlled dilations of one of the families
of simplices.
\end{rem}

Now, from the fact that $\sigma_1 \pitchfork_{\delta} \sigma_2$ it follows that $\sigma_3 \cap \sigma_4 \neq
\emptyset \Leftrightarrow \bar{\sigma}_3 \cap \bar{\sigma}_4 \neq \emptyset$ (see \cite{cms}, p. 436). Therefore,
the function $Vol_{eucl}(\gamma_0)$ attains a positive minimum, as a positive, continuous function defined on the
compact set $\bar{\sigma}_3 \cap \bar{\sigma}_4$, thus proving the first assertion of the proposition.

Let be $\gamma$ be a $q$-dimensional cell, and let $\beta$ be a face of $\partial\gamma$. Then:

\begin{equation}
Vol_{eucl}(\beta) \leq d_1^{p}\;.
\end{equation}

Choose $p \in \gamma$, such that $\rho = dist(p,\partial\gamma) = \max\{dist(r,\partial\gamma)\,|\,r \in
\gamma\}$. Then, if $\beta = \beta^j$ denotes a $j$-dimensional face of $\partial\gamma$, we have that:

\begin{equation}
\gamma \subset \bigcup_{\beta^j \subset \partial\gamma}\!\!N_{\rho}(\beta^j);
\end{equation}

where: $N_{\rho}(\beta^j) = \{r\,|\, dist(r,\beta^j) \subseteq \rho\}$. But:

\begin{equation}
Vol_{eucl}(\beta^j \cap \gamma) \leq c\cdot\rho^{q-j}\cdot Vol_{eucl}(\beta^j),
\end{equation}

for some $c' = c'(q)$.
\\ Moreover, the number of faces $\sigma_3 \cap \sigma_4$ of $\gamma$ is at most
$2^{dim\,\sigma_1+dim\,\sigma_2+2}$, where $\sigma_1, \sigma_2$ are as in Remark 3.4. and  $dim\,\sigma_3 \leq
dim\,\sigma_1,\;  dim\,\sigma_4 \leq dim\,\sigma_2$.
\\Thus (3.4) in conjunction with (3.6) imply that there exists $c_1 = c_1(m,\varphi_0,\delta)$, such that:

\begin{equation}
c_1d_1^q \leq \sum_{j=0}^{q-1}\rho^{q-j}d_1^j.\,
\end{equation}

and (3.2) follows from this last inequality.

The last inequality follows from (3.2) and (3.3) by induction.

\end{proof}

Next we show that given two fat Euclidian triangulations that intersect $\delta$-transver-sally, then one can
infinitesimally move any given point of one of the triangulations s.t. the resulting intersection will be
$\delta'$-transversal, where $\delta'$ depends only on  $\delta$, the common fatness of the given triangulations,
and on the displacement length. More precisely one can show that the following results holds:

\begin{prop}
Let $K_1, K_2 \subset \mathbb{R}^n$ be $n$-dimensional simplicial complexes, of common fatness $\varphi_0$ and
$d_1 = diam\,\sigma_1 \leq d_2 = diam\,\sigma_2$. Let $v_0 \in K_1$ be a $0$-dimensional simplex of $K_1$.
Consider the complex $K^*_1$ obtained by replacing $v_0$ by $v^*_0 \in R^n$ and keeping fixed the rest of the
$0$-dimensional vertices of $K_1$ fixed. Consider also $L_2 < K_2,\; L_2 = \{\sigma \in K_2\,|\, \sigma \cap
B(v_0,2d_1) \neq \emptyset\}$.
\\ Then,  if there exists $k$ s.t. all the $k$-simplices $\tau \subset \partial St(v_0)$ are $\delta$-transversal
to $L_2$, there exist $\varphi_0, \delta, \varepsilon > 0$, $\delta^\ast =
\delta^\ast(\varphi_0,\delta,\varepsilon)$ and there exists $v^*_0$ s.t. $dist(v_0,v^*_0) < \varepsilon \cdot d_1$
s.t.

\begin{equation}
\tau^* \pitchfork_{\delta^{\ast}} L_2\,; \;\forall\, \tau \subset St(v^*_0)\, \backslash\, \partial St(v^*_0),\;
dim\, \tau^* = k+1.
\end{equation}

\end{prop}

\begin{proof}
Let $N(r) = |\{ \sigma \in K_1\,|\, \sigma \subset B_r(v_0)\}|$. Then there exists a constant $c_n$ s.t. $N(r)
\leq \frac{c_n}{\varphi_0}(\frac{\varepsilon}{d_1})^n$. It follows that the set $St(v_0)$ is compact, since there
are at most $\frac{c_n}{\varphi_0}$ possible edge lengths, which can take values in the interval
$[d_1\varphi_0,d_1]$.\footnote{\,i.e. the number of possible combinatorial structures on $St(v_0)$ depends only on
$\varphi_0$.} Therefore if a $D^\ast$ satisfying (3.8) exist, it depends only on $\varphi_0,\; \delta$ and
$\varepsilon$ (and not on $K_1, K_2$).
\\ Let $\sigma_1, \ldots , \sigma_{l_1}$ and $\tau_1, \ldots, \tau_{l_2}$ be orderings of the simplices of $L_2$ and of the $k$-simplices of $\partial St(v_0)$, respectively.
Then, by \cite{cms}, Lemma 7.4, there exists $\varepsilon_{1,1}$ and exists $v_{1,1}, \; d(v_{1,1},v_0) =
\varepsilon_{1,1}$, such that the hyperplane $\Pi(v_{1,1},\tau_1)$ determined by $v_1$ and by $\tau_1$ is
transversal to $\sigma_1$. By replacing $\tau_1$ by $\tau_2$ and $v_0$ by $v_{1,1}$ we obtain $v_{1,2}$ and
$\varepsilon_{1,2}$ s.t. $\Pi(v_{1,2},\tau_1) \pitchfork \sigma_2$. (See Fig. 3.)

\begin{figure}[h]
\begin{center}
\includegraphics[scale=0.35]{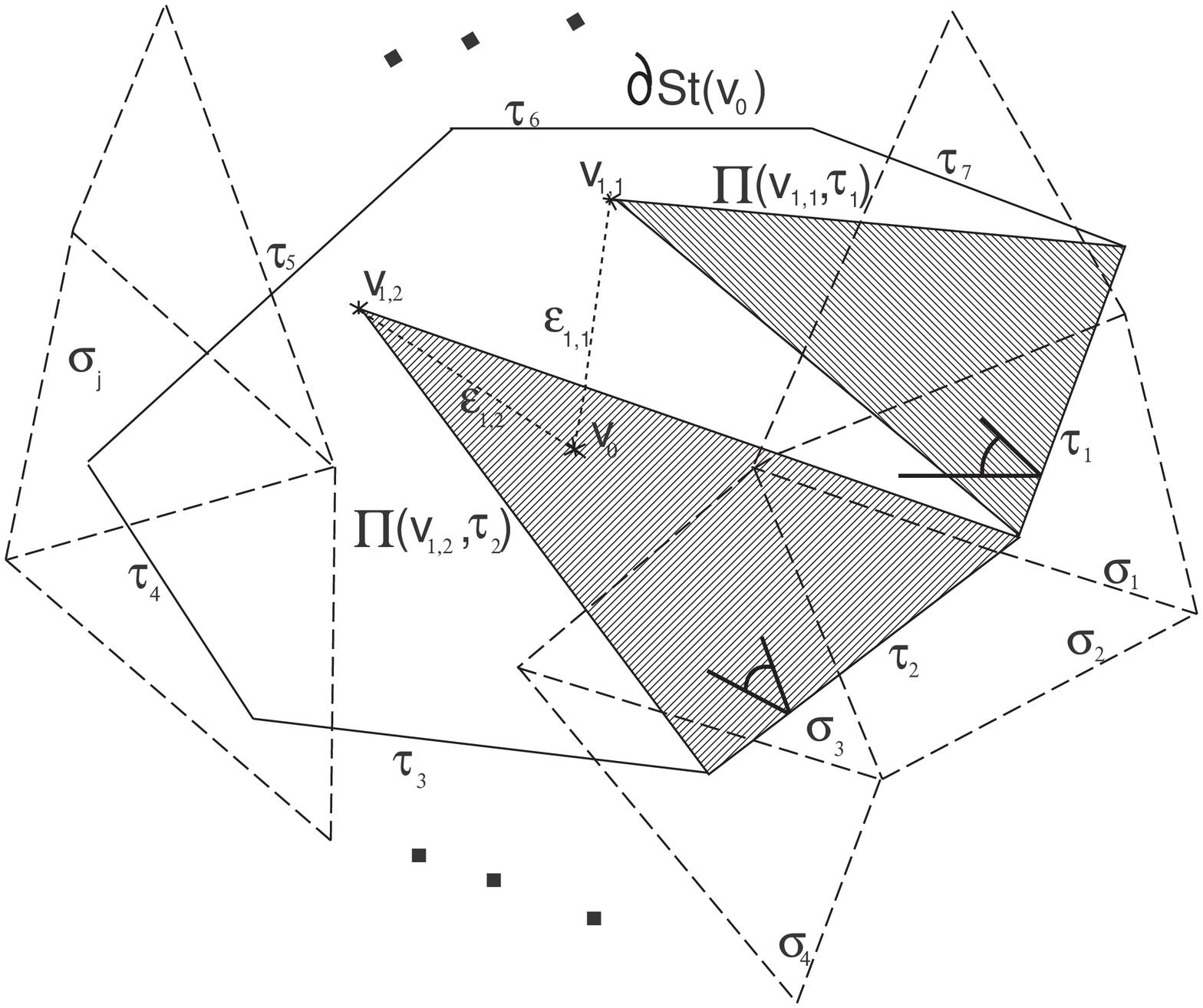}
\end{center}
\caption{}
\end{figure}

 Moreover, by choosing
$\varepsilon_{1,2}$ sufficiently small, one can ensure that $\Pi(v_{1,2},\tau_1) \pitchfork \sigma_1$, also.
Repeating the process for $\tau_3, ..., \tau_{l_2}$, one determines a point  $v_{1,l_2}$ such that
$\Pi(v_{1,l_2},\tau_j) \pitchfork L_2\; j = 3,\ldots,l_2$. In the same manner and by choosing at each stage an
$\varepsilon_{i,j}$ small enough, one finds points $v_{i,j}$ s.t. $\Pi(v_{i,j},\tau_i) \pitchfork L_2,\; i = 1,
\ldots ,l_1\,,\; j = 1,\dots,l_2$. Then $v^*_0 = v_{l_1,l_2}$ satisfies: $\Pi(v_0^\ast,\tau_{j}) \pitchfork L_2$,
$j = 1,\ldots,,l_2$.

\end{proof}

We are now prepared to prove the main result of this section namely:

\begin{thm}
Let $\mathcal{T}_1,\mathcal{T}_2$ be two fat triangulations of open sets $U_1,U_2 \subset M^n,\,$ $U_1 \cap U_2
\neq \emptyset$ having common fatness $\geq \varphi_0$, and such that  $\mathcal{T}_1 \cap \mathcal{T}_2 \neq
\emptyset$. Then there exist fat triangulations $\widetilde{\mathcal{T}'_1},\mathcal{T}'_2$ and there exist open
sets $U \subset  U_1 \cap U_2 \subset V$, such that
\begin{enumerate}
\item $({\mathcal{T}'_1} \cap \mathcal{T}'_2) \cap (U_i \setminus V) = \mathcal{T}_i\,,\; i=1,2\,;$
\item $({\mathcal{T}'_1} \cap \mathcal{T}'_2) \cap U = \mathcal{T};$
\\\hspace*{-1.3cm} where
\item $\mathcal{T}$ is a fat triangulation of $U$.
\end{enumerate}
\end{thm}

\begin{proof}
Let $K_1,\,K_2$ denote the underlying complexes of $\mathcal{T}_1,\mathcal{T}_2$, respectively. By considerations
similar to those of Proposition 3.5. it follows that given $\varphi_0
> 0$, there exists $d(\varphi_0) > 0$ such that given a $k$-dimensional simplex $\sigma \subset \mathbb{R}^n$,
$diam(\sigma) = d_1$ has fatness $\varphi_0$, than translating each vertex of $\sigma$ by a distance
$d(\varphi_0)\cdot d_1$ renders a simplex of fatness $\geq \varphi_0/2$\,. Also, it follows that given
$\varphi_0,\delta
> 0$, exists $\delta(\varphi_0,\delta)$ satisfying the following condition: if every vertex $u \in \sigma \subset
K$ is replaced by a vertex $u'$ s.t. $dist(u,u') \leq \delta(\varphi_0,\delta) \cdot d_1$, then the resulting
simplex $\sigma'$ is $\pitchfork_{\delta/2}$-transversal to $K$; for any $n$-dimensional simplicial complex $K$ of
fatness $\varphi_0$ and such that $diam\,\sigma = d_2 \geq d_1$ .
\\ Let $v_0 \in U_1 \cap U_2$. Define the following subcomplexes of $K_1,\,K_2$, respectively:
\[L_2 = \{\bar{\sigma} \subset K_2 \,|\, \bar{\sigma} \subset B_\varepsilon(v_0),\, d_1 \leq dist(\bar{\sigma},\partial B_\varepsilon(v_0)) \leq d_2\}\]
\[M_2 = \{\sigma \subset K_2 \,|\, \bar{\sigma} \subset \bar{\tau} \subset B_\varepsilon(v_0),\, dim\,\tau = n,\,\bar{\sigma} \cap L_2 \neq \emptyset\}\]
\[\hspace*{-3.6cm}L_1 = \{\bar{\sigma} \subset K_1 \,|\, dist(\bar{\sigma},L_2) \leq d_2\}\]
\[M_1 = \{\sigma \subset K_1 \,|\, \bar{\sigma} \subset \tau \subset B_\varepsilon(v_0),\, dim\,\tau = n,\, \tau \cap L_1 \neq \emptyset\}\]
(See Fig. 4.)
\begin{figure}[h]
\begin{center}
\includegraphics[scale=0.35]{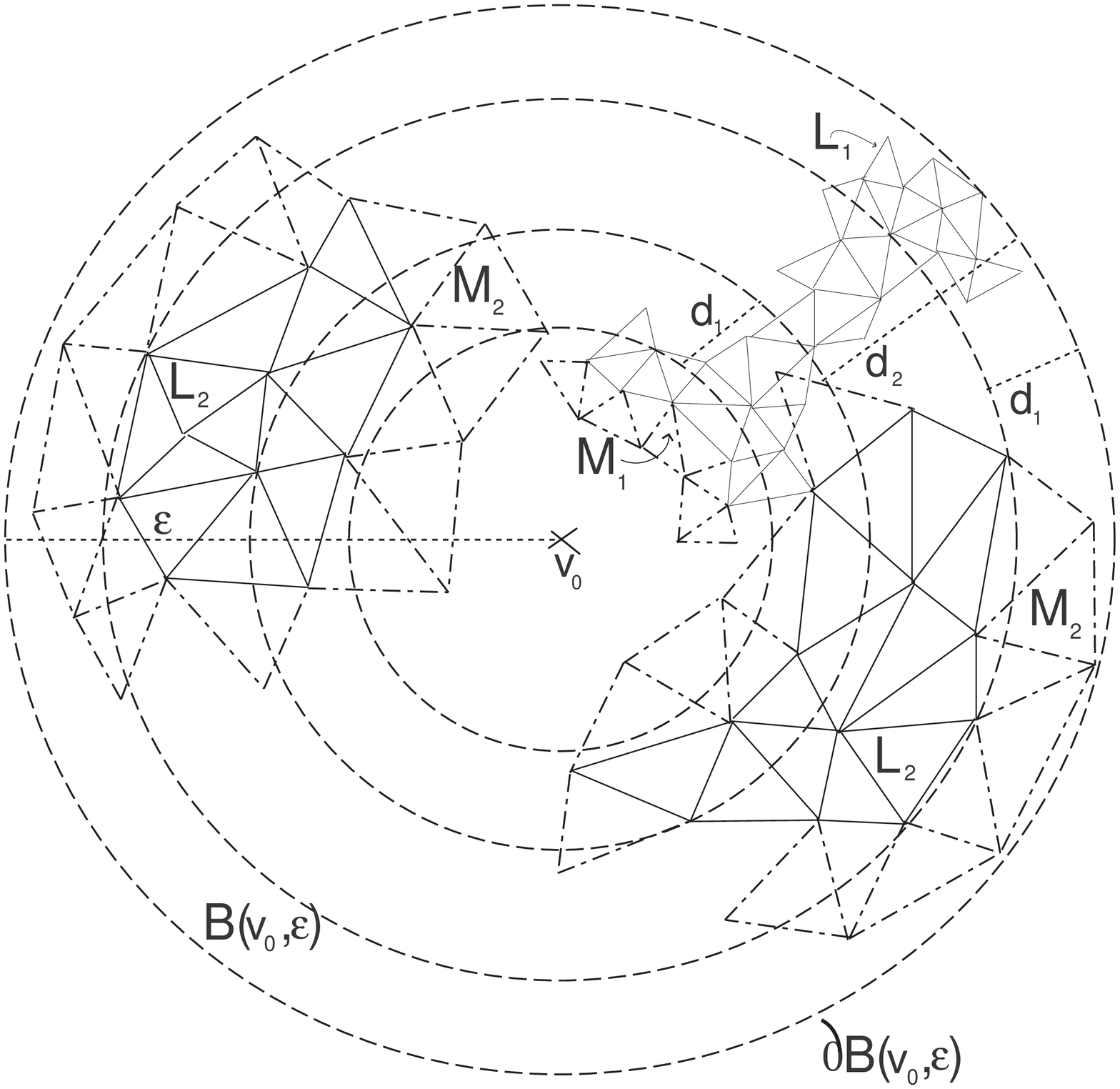}
\end{center}
\caption{}
\end{figure}

Consider an ordering $v_1, \ldots, v_p$ of the vertices of $L_1$. It follows from Proposition \nolinebreak[4]3.5.
that, if all the vertices of $L_1$ are moved by at most $t_0$, where

\begin{equation}
t_0 = \frac{d_1}{n}\min{\big\{\frac{1}{2},d(\varphi_0)\big\}}\,,
\end{equation}

then there exists

\begin{equation}
\delta^\ast_0 = \delta_0^\ast(\varphi_0,1,\frac{t_0}{d_1})\,,
\end{equation}

such that

\begin{equation}
\mathcal{S}^0(L_{1,0}) \pitchfork_{\delta_{0}^\ast} K_2\,,
\end{equation}

where $\mathcal{S}^0(L_{1,0})$ denotes the $0$-skeleton of $L_{1,0}$.

Now define inductively

\begin{equation}
t_i =\frac{d_1}{n}\min{\bigg\{\frac{1}{2},d(\varphi_0),\delta\Big(\frac{\varphi_0}{2},\frac{\delta^\ast_0}{2}\Big)
,\ldots,\delta\Big(\frac{\varphi_0}{2},\frac{\delta^\ast_{i-1}}{2}\Big) \bigg\}}\,,
\end{equation}

where

\begin{equation}
\delta^\ast_{i} = \delta^\ast_{i}\Big(\varphi_0,\delta^\ast_{i-1},\frac{t_i}{d_1}\Big)\,;\; i=1,\ldots,n-1.
\end{equation}

Then $t_0 \geq t_1 \geq \ldots \geq t_{n-1}$.
\\ Moving  each and every vertex of $L_1$ by a distance $\leq t_i,\;i=1,\ldots,n$, renders complexes
$L_{1,1},\ldots,L_{1,n-1}$ s.t.
\begin{enumerate}
\item $L_{1,i}\cap\big(B_\varepsilon(v_0) \setminus \, M_1 \big) \equiv L_1$,
\item $L_{1,i}$ are $\varphi_0$-fat; $i=0,\ldots,n-1$.
\end{enumerate}

By inductively applying Proposition 3.5. it follows that

\begin{equation}
\mathcal{S}^i(L_{1,i}) \pitchfork_{\delta_i^\ast} K_2\,,
\end{equation}

Where $\mathcal{S}^i(L_{1,i})$ denotes the $i$-skeleton of $L_{1,i}$\,.
\\Moreover,

\begin{equation}
\mathcal{S}^i(L_{1,j}) \pitchfork_{\delta_i^\ast/2} K_2\,, \forall j > i.
\end{equation}

It follows that

\begin{equation}
L_{1,n-1} \pitchfork_{\delta^\bigstar} K_2\,,
\end{equation}

where

\begin{equation}
\delta^\bigstar = \frac{1}{2}\min\{\delta_0^\ast,\ldots,\delta_i^\ast\}\,.
\end{equation}

By Proposition 3.2. the barycentric subdivision of $L_{1,n-1} \cap L_2$ is fat. We extend it to a fat subdivision
of $M_2$ in the following manner: given a simplex $\sigma \subset  M_2 \setminus L_2$, subdivide $\sigma$ by
constructing all the simplices with vertices $v_i$, where $v_i$ is either the vertex of a simplex $\sigma \subset
M_2 \setminus L_2,\; \bar{\sigma_i} \cap L_2 \neq \emptyset$, or it is a vertex of a closed simplex $\bar{\sigma}$
of the barycentric subdivision of $L_{1,n-1} \cap L_2$\,, such that $\bar{\sigma} \subset \partial L_2 \cap
M_2\,,\; i = 1,\ldots,k_0$. The triangulation $\widetilde{K}_2$ thus obtained is a fat extension of $\overline{K_2
\setminus\,M_2 }$. 

In an analogous  manner one constructs a similar fat extension $\widetilde{K}_1$ of $\overline{K_1
\setminus\,M_2}$.

\end{proof}

Now let $\mathcal{T}_1$, $\mathcal{T}_2$ be the triangulations of $\partial M^n \times [0,1)$ and $int\,M^n$,
respectively, given by Theorem 2.9.\,. Then the local fat triangulation obtained in Theorem 3.6. extends globally
to a fat triangulation of $\mathcal{T}_1 \cap \mathcal{T}_2$, by applying Lemma 10.2 and Theorem  10.4 of
\cite{mun}. This concludes the
\[\hspace{8.5cm} \fbox{Proof of Theorem 1.1}\]

\section{The Existence of Quasimeromorphic Mappings}

\subsection{Quasimeromorphic Mappings}

\begin{defn}
Let $D \subseteq \mathbb{R}^n$ be a domain; $n \geq 2$, and let $f: D \rightarrow \mathbb{R}^m$.
\\ $f$ is called $ACL$ ({\it absolutely continuous on lines}) iff:
\\ (i) $f$ is continuous
\\ (ii) for any $n$-interval $Q = \overline{Q} = \{a_i \leq x_i \leq b_i\,|\, i=1,\ldots,n\}$, $f$ is absolutely
continuous on almost every line segment in $Q$, parallel to the coordinate axes.
\end{defn}

\begin{lem} [\cite{v}, 26.4]
If $f:D \subseteq \mathbb{R}^n \rightarrow  \mathbb{R}^m$ is $ACL$, then $f$ admits partial derivatives almost
everywhere.
\end{lem}

The result above justifies the following Definition:

\begin{defn}
$f:D \subseteq \mathbb{R}^n  \rightarrow  \mathbb{R}^m$ is $ACL^p$ iff its derivatives are locally $L^p$
integrable, $p \geq 1$.
\end{defn}

\begin{defn}
Let $D \subseteq \mathbb{R}^n$ be a domain; $n \geq 2$ and let $f: D \rightarrow \mathbb{R}^m$ be a continuous
mapping. $f$ is called
\begin{enumerate}
\item {\it quasiregular} iff (i) $f$ is $ACL^n$ and
\\ \hspace*{2.3cm}(ii) $\exists K \geq 1$ s.t.
\begin{equation}
|f'(x)| \leq KJ_f(x)\; a.e.
\end{equation}
\\ where $f'(x)$ denotes the formal derivative of $f$ at $x$, $|f'(x)| = \sup \raisebox{-0.2cm}{\mbox{\hspace{-0.7cm}\tiny$|h|=1$}}|f'(x)h|$, and where $J_f(x) = detf'(x)$.
\\ The smallest $K$ that satisfies (4.1) is called the {\it outer dilatation} of $f$.
\item {\it quasiconformal} iff $f:D \rightarrow f(D)$ is a local homeomorphism.
\item {\it quasimeromorphic} iff $f:D \rightarrow \widehat{\mathbb{R}^n}$, $\widehat{\mathbb{R}^n} = \mathbb{R}^n \bigcup
\,\{\infty\}$ is quasiregular, where the condition of quasiregularity at $f^{-1}(\infty)$ can be checked by
conjugation with auxiliary M\"{o}bius transformations.
\end{enumerate}
\end{defn}

\begin{rem}
One can extend the definitions above to oriented $\mathcal{C}^\infty$ Riemannian $n$-manifolds by using coordinate
charts.
\end{rem}

\subsection{Alexander's Trick}
The technical ingredient in Alexander's trick is the following Lemma:
\begin{lem} (\cite{ms1}, \cite{pe})
Let $\mathcal{T}$ be a fat triangulation of $M^n \subset \mathbb{R}^n$, and let $\tau,\sigma \in
\nolinebreak[4]\mathcal{T},\; \tau = (p_1,\dots,p_n), \,\sigma = (q_1,\dots,q_n)$; and denote $|\tau| = \tau \cup
int\,\tau$.
\\ Then there exists a sense-preserving homeomorphism $h = h_{\tau}: |\tau| \rightarrow \widehat{\mathbb{R}^{n}}$
s.t.
\begin{enumerate}
\item $h(|\tau|) = |\sigma|$, \,if\; $det(p_1,\dots,p_n) > 0$
\\ and
\\ $h(|\tau|) = \widehat{\mathbb{R}^{n}} \setminus| \sigma|$, \,if\; $det(p_1,\dots,p_n) < 0$.
\item $h(p_i) = q_i, \; i=1,\ldots,n.$
\item $h|_{\partial|\sigma|}$ is a $PL$-homeomorphism.
\item $h|_{int|\sigma|}$ is quasiconformal.
\end{enumerate}
\end{lem}

\begin{prf}
Let $\tau_0 = (p_{0,1},\dots,p_{0,n})$ denote the equilateral $n$-simplex inscribed in the unit sphere
$\mathbb{S}^{n-1}$. The {\it radial linear stretching} $\varphi: \tau \rightarrow \overline{B^n}$ is onto and
bi-lipschitz (see \cite{ms2}). Moreover, by a result of Gehring and V\"{a}isal\"{a}, $\varphi$ is also
quasicomformal (see \cite{v}). We can extend $\varphi$ to $\widehat{\mathbb{R}^n}$ by defining $\varphi(\infty) =
\infty$.
\end{prf}

The Existence Theorem  of quasimeromorphic mappings now follows immediately:

\begin{thm}
Let $M^n$ be a $C^\infty$ Riemannian manifold with or without boundary. Then there exists a quasimeromorphic
mapping $f:M^n \rightarrow \widehat{\mathbb{R}^n}$.
\end{thm}

\begin{prf}
Let $f:M ^n\rightarrow \widehat{\mathbb{R}^n}$ be defined by: $f|_{|\sigma|} = h_{\sigma}$, where $h$ is the
homeomorphism constructed in the Lemma above. Then $f$ is a local homeomorphism on the $(n-1)$-skeleton of
$\mathcal{T}$ too, while its branching set $B_f$ is the $(n-2)$-skeleton of $\mathcal{T}$. By its construction $f$
is quasiregular. Moreover, given the uniform fatness of the triangulation $\mathcal{T}$, the dilatation of $f$
depends only on the dimension \nolinebreak[4]$n$.
\end{prf}

\section{Generalizations and Further Research}

We succinctly present some immediate generalizations of Theorems 1.1. and 2.9.

\subsection{Smoothings}
Theorem  1.1. was restricted to  $\mathcal{C}^{\infty}$ manifolds because the  triangulation $\mathcal{T}_2$ of
$int\,M^n$ was obtained by applying Peltonen's Theorem; so our overall argument is valid only for
$\mathcal{C}^{\infty}$ manifolds. But  the class of any $n$-manifold may be elevated up to $\mathcal{C}^{\infty}$
(see \cite{mun}, Theorems 4.8 and 5.13 ), so we can apply the methods of \cite{pe} on the smoothed
$\mathcal{C}^{\infty}$ manifold, and then project the fat triangulation received to the original structure. Since
in the smoothing process we employed only $\delta$-approximations that are, by Lemma 8.7. \cite{mun},
$\alpha$-approximations too, we will obtain a fat triangulation, as desired. We can thus formulate the following
Corollary:

\begin{cor}
Let $M^{n}$ be an $n$-dimensional $\mathcal{C}^{r}\,, 1 \leq r \leq \infty$ manifold with boundary. Then any fat
triangulation of $\partial M^{n}$ can be extended to a fat triangulation of $M^{n}$.
\end{cor}

Moreover, every $PL$ manifold of dimension $n \leq 4$ admits a (unique, for $n \leq 3$) smoothing (see
\cite{mun1}, \cite{mun}, \cite{th2}), and every topological manifold of dimension $n \leq 3$ admits a $PL$
structure (cf. \cite{moi}, \cite{th2}). Therefore we can can start with a $PL$ manifold (or even just a
topological one in dimensions $2$ and $3$) and smooth it, thus receiving

\begin{cor}
Let $M^{n}$ be an $n$-dimensional, $n \leq 4$ (resp. $n \leq 3$), $PL$ (resp. topological) manifold with boundary.
Then any fat triangulation of $\partial M^{n}$ can be extended to a fat triangulation of $M^{n}$.
\end{cor}

Using again Alexander's Trick renders the following result:

\begin{cor}
Let $M^n$ be an $n$-dimensional manifold ($n \geq 2$),  with or without boundary. Then in the following cases
there exists a quasimeromorphic mapping $f:M^n \rightarrow \widehat{\mathbb{R}^n}$:
\begin{enumerate}
\item $M^n$ is of class $\mathcal{C}^{r}\,, 1 \leq r \leq \infty\,, \; \forall n \geq 2$;
\item $M^n$ is a $PL$ manifold and $n \leq 4$;
\item $M^n$ is a topological manifold and $n \leq 3$.
\end{enumerate}
\end{cor}

\begin{rem}
 It may be that during the smoothing process the dilatation will increase in an unbounded fashion, thus
rendering impossible the proof of existence of fat triangulations quasimeromorphic mappings. However, the
dilatation may increase only when we {\it linearize} the {\it tangent cone} at {\it cone points}. Fortunately, the
nature of linearization process is such that, when the cone angles are bounded from below, then the dilatations
will be bounded from above\footnote{\,For some details regarding linearizations,  see \cite{th2}.}.
\end{rem}

\subsection{Kleinian Groups}
Since the construction of fat triangulations was motivated mainly by the study of $G$-automorphic quasimeromorphic
mappings with respect to a {Kleinian group} $G$, i.e. a discontinuous group of orientation preserving isometries
of $\mathbb{H}^n$, it is natural to employ Theorems 1.1. and 2.9. to prove the following result, that represents a
generalization of a result of Tukia (\cite{tu}):

\begin{thm} Let $G$ be a Kleinian group with torsion acting upon $\mathbb{H}^n, \, n \geq 3$.\footnote{\, The case $n =2$ being trivial, since in this case $\mathbb{H}^n/G$ is always a manifold.} \\ If the elliptic elements (i.e. torsion elements) of $G$ have
uniformly bounded orders,
\\then there exists a non constant $G$-automorphic quasimeromorphic mapping
\\$f: \mathbb{H}^n \rightarrow \widehat{\mathbb{R}^n}$, i.e. such that
 \begin{equation} f(g(x))= f(x)\,,\;\forall x \in \mathbb{H}^n\,, \forall g \in G\,.
 \end{equation}
\end{thm}

While for full details we refer the reader to \cite{s2} and -- for a different fattening method (albeit in
dimension 3 only), to \cite{s1} -- we bring here the following

\begin{skprf} By Lemma 4.6. it suffices to produce a fat $G$-invariant triangulation of $\mathbb{H}^n$.
The {\em singular locus} $\mathcal{L}$ of $\mathbb{H}^n/G$ is the image, under the natural projection
$\pi:\mathbb{H}^n \rightarrow \mathbb{H}^n/G$\,, of the union $\mathcal{A} =
\bigcup\raisebox{-0.7em}{\hspace{-0.5cm}\tiny $i\in \mathbb{N}$}A_{f_i}$ of the elliptic axes of $G$\footnote{\,
$\mathcal{A} = \{A_{f_i}\}_{_i}$ is a countable set, by the discreteness of $G$.}. For each elliptic axes
$A_{f_i}$ it is possible to choose a {\em collar} $N_i$ and triangulate it in an $f_i$-invariant manner. Denote by
$\mathcal{T}_i$ the  $f_i$-invariant triangulation of $N_i$.  Put $\mathcal{N} =
\bigcup\raisebox{-0.7em}{\hspace{-0.5cm}\tiny $i\in \mathbb{N}$}N_i$. Then $M_e = (\mathbb{H}^n \backslash
\,\mathcal{N})/\,G $ is a manifold with boundary. Then $\partial M_e =
\bigcup\raisebox{-0.7em}{\hspace{-0.5cm}\tiny $i\in \mathbb{N}$}\partial N_i$ has the triangulation induced by
that of $\mathcal{N}$\footnote{\, Some more care is needed in producing the fat triangulation of $\mathcal{N}$ if
there exist $i \neq j$ s.t. $A_{f_i} \cap A_{f_j} \neq \emptyset$.} and, since the orders of the elliptic elements
are bounded from above, the induced triangulation will be fat. By Theorem 1.1.\,, this triangulation can be
extended to a fat triangulation $\mathcal{T}$ of $\mathbb{H}^n/G$\,. Then $\pi^{-1}(\mathcal{T}) \,\cup\,
\bigcup\raisebox{-0.7em}{\hspace{-0.5cm}\tiny $i\in \mathbb{N}$} \mathcal{T}_i$ will represent the desired fat
$G$-invariant triangulation.
\end{skprf}

\begin{rem}
It seems feasible to adapt this argument for any {\it geometric orbifold} with tame singular locus (at least in
dimension 3).
\end{rem}

\subsection{Lipschitz Manifolds}
The existence of triangulations for {\it Lipschitz manifolds} was already stated by Cairns (\cite{ca3}), yet it
was never proved in full detail. Also it seems possible relax the smoothness condition even further, as to include
{\it quasisymmetric} manifolds, i.e. manifolds of which local charts are given by {\it quasisymmetric} mappings,
where:

\begin{defn}
An embedding $f:\mathbb{R}^m  \rightarrow \mathbb{R}^n$ is called {\it quasisymmetric} iff here exists a
homeomorphism $\eta: [0,\infty) \rightarrow [0,\infty)$ such that for all $x,a,b \in \mathbb{R}^m $ and for all $t
\in [0,\infty)$ the following holds:

\begin{equation}
dist(x,a) \leq t\cdot dist(x,b) \Longrightarrow dist(f(x),f(a)) \leq \eta(t)\cdot dist(f(x),f(b))\,.
\end{equation}

\end{defn}


\end{document}